\documentstyle{amsppt}
\baselineskip18pt
\magnification=\magstep1
\pagewidth{4.5in}
\pageheight{7.0in}
\hyphenation{co-deter-min-ant co-deter-min-ants pa-ra-met-rised
pre-print fel-low-ship Cox-et-er semi-stand-ard}
\def\leaderfill{\leaders\hbox to 1em{\hss.\hss}\hfill}

\def\idest{i.e.\ }

\def\l{{\lambda}}

\def\rad{\text{\rm rad}}

\def\boxit#1{\vbox{\hrule\hbox{\vrule \kern3pt
\vbox{\kern3pt\hbox{#1}\kern3pt}\kern3pt\vrule}\hrule}}
\def\rabbit{\vbox{\hbox{\kern0pt
\vbox{\kern0pt{\hbox{---}}\kern3.5pt}}}}

\def\tableau#1{
        \hbox {
                \hskip -10pt plus0pt minus0pt
                \raise\baselineskip\hbox{
                \offinterlineskip
                \hbox{#1}}
                \hskip0.25em
        }
}

\def\tabCol#1{
\hbox{\vtop{\hrule
\halign{\strut\vrule\hskip0.5em##\hskip0.5em\hfill\vrule\cr\lower0pt
\hbox\bgroup$#1$\egroup \cr}
\hrule
} } \hskip -10.5pt plus0pt minus0pt}

\def\CR{
        $\egroup\cr
        \noalign{\hrule}
        \lower0pt\hbox\bgroup$
}



\def\xcarti{{\bf 1}}

\def\xdjb{{\bf 2}}

\def\xdub{{\bf 3}}

\def\xduf{{\bf 5}}

\def\xduh{{\bf 4}}
\def\xdr{{\bf 6}}
\def\xeg{{\bf 7}}

\def\xgl{{\bf 8}}

\def\xjagb{{\bf 9}}

\def\xrmge{{\bf 11}}

\def\xthesis{{\bf 10}}

\def\xkx{{\bf 12}}

\def\xlusc{{\bf 13}}

\topmatter
\title Completions of cellular algebras
\endtitle

\author R.M. Green\endauthor
\affil 
Department of Mathematics and Statistics\\ Lancaster University\\
Lancaster LA1 4YF\\ England\\
{\it  E-mail:} r.m.green\@lancaster.ac.uk
\endaffil

\abstract
We introduce procellular algebras, so called because they are
inverse limits of finite dimensional cellular algebras as defined by
Graham and Lehrer.  A procellular algebra is defined as a certain 
completion of an infinite dimensional cellular algebra whose cell
datum is of ``profinite type''.  We show how these notions overcome
some known obstructions to the theory of cellular algebras in infinite 
dimensions.
\endabstract

\thanks
The author was supported in part by an E.P.S.R.C. postdoctoral
research as\-sist\-ant\-ship.
\endthanks
\endtopmatter

\centerline{\bf To appear in Communications in Algebra}

\head 0. Introduction \endhead

Cellular algebras, which were introduced by Graham and Lehrer in
[\xgl], are a class of associative algebras defined in terms
of a ``cell datum'' satisfying certain axioms.  Several interesting
classes of finite dimensional algebras, including certain Hecke
algebras and Brauer algebras, can be described in this way.  

One of the strengths
of the theory of cellular algebras is that it provides a complete list
of absolutely irreducible modules for the algebra 
over a field.  Unfortunately, this property is
no longer true if one allows an algebra satisfying the cellular axioms
to be infinite dimensional, and there is a simple counterexample
to show this.  

The goal of this paper is to give a natural definition of cellularity for a
general algebra which not only agrees with the original definition if
the algebra is finite dimensional, but which also produces results analogous to
those in the finite dimensional case if the algebra is infinite dimensional.

The material in \S1 is expository, giving a summary of some of the
important results in the theory of cellular algebras, and outlining the
obstructions to the theory in infinite dimensions.  The main problem
is that although the definition of ``cellular'' makes sense for
infinite as well as finite dimensional algebras, the results of [\xgl]
assume tacitly that the algebra is finite dimensional.  If the algebra
is infinite dimensional, the results can fail badly even for simple
examples.

In \S2, we introduce the notion of a cell datum of ``profinite
type''.  For such a cell datum, the associated partially ordered set
is allowed to be infinite, provided that it satisfies a certain
finitary condition.  If the condition is satisfied, the infinite
dimensional algebra of profinite type may be completed with respect to
a natural topology, which depends on the cell datum, resulting in a
``procellular'' algebra.  The desired classification for absolutely
irreducible modules can then be obtained by considering ``smooth'' 
representations of the completed algebra.

In \S3, we look at a nontrivial natural example of such a procellular
algebra.  This example arises as a completion of Lusztig's algebra
$\dot U$, which
appears in [\xlusc], or alternatively as a projective limit of
$q$-Schur algebras.  We also discuss briefly another very interesting 
natural example of a procellular algebra which arises from a certain
completion of the affine Temperley--Lieb algebra.

\vskip 20pt

\head 1. Cellular algebras \endhead

\subhead 1.1 Definitions \endsubhead

We start by reviewing some of the basic properties of cellular algebras.

\proclaim{Definition 1.1.1 ([\xgl, Definition 1.1])}  Let $R$ be a commutative
ring with identity.  A {\it cellular algebra} over $R$ is an associative unital
algebra, $A$, together with a cell datum $(\Lambda, M, C, *)$ where

\item {\rm 1.}
{$\Lambda$ is a partially ordered set (or ``poset'' for short), ordered
by $<$.  For each $\l \in \Lambda$, $M(\l)$ is a finite set
(the set of ``tableaux'' of type $\l$) such that $$
C : \coprod_{\l \in \Lambda} \left( M(\l) \times M(\l) \right) \rightarrow A
$$ is injective with image an $R$-basis of $A$.}
\item {\rm 2.}
{If $\l \in \Lambda$ and $S, T \in M(\l)$, we write $C(S, T) = C_{S, T}^{\l}
\in A$.  Then $*$ is an $R$-linear involutory anti-auto\-morph\-ism 
of $A$ such that
$(C_{S, T}^{\l})^* = C_{T, S}^{\l}$.}
\item {\rm 3.}
{If $\l \in \Lambda$ and $S, T \in M(\l)$ then for all $a \in A$ we have $$
a C_{S, T}^{\l} \equiv \sum_{S' \in M(\l)} r_a (S', S) C_{S', T}^{\l}
\mod A(< \l),
$$ where  $r_a (S', S) \in R$ is independent of $T$ and $A(< \l)$ is the
$R$-submodule of $A$ generated by the set $$
\{ C_{S'', T''}^{\mu} : \mu < \l, S'' \in M(\mu), T'' \in M(\mu) \}
.$$}
\endproclaim

We use Definition 1.1.1 only when $\Lambda$ is a finite set, although
such an assumption is omitted in [\xgl].

There are many interesting examples of finite dimensional algebras
which satisfy these axioms, some of which are described in detail in
[\xgl].  These include the Brauer algebra,
the Temperley--Lieb algebra, the Hecke algebra of type $A$ (using the
Kazhdan--Lusztig basis as the cell basis) and Jones' annular algebra.

The theory of cellular algebras may be applied to find, among other
things, criteria for semisimplicity and
a complete set of absolutely irreducible modules over a
field.  In this paper, we will mostly be concerned with the latter
application, and to understand it, we recall the notion of cell
modules.

\proclaim{Definition 1.1.2 ([\xgl, Definition 2.1]}
For each $\l \in \Lambda$ define the left $A$-module $W(\l)$ as
follows: $W(\l)$ is a free $R$-module with basis $\{C_S: S \in
M(\l)\}$ and $A$-action defined by $$
a C_S = \sum_{S' \in M(\l)} r_a(S', S)C_{S'} \quad (a \in A, S \in M(\l))
,$$ where the $r_a$ are the structure constants of Definition 1.1.1.
This is called the cell representation of $A$ corresponding to $\l \in
\Lambda$.
\endproclaim

Associated with each cell representation $W(\l)$ is a bilinear form
$\phi_{\l}$ which is of key importance in classifying the irreducible
modules.

\proclaim{Lemma 1.1.3}
Let $\l \in \Lambda$.  Then for any elements $S_1, S_2, T_1, T_2 \in
M(\l)$, we have $$
C_{S_1, T_1}^{\l} C_{S_2, T_2}^{\l} \equiv \phi_1(T_1, S_2) C_{S_1,
T_2}^{\l} \mod A(< \l)
$$ where $\phi_1 \in R$ is independent of $T_2$ and $S_1$.
\endproclaim

\demo{Proof}
This follows from [\xgl, Lemma 1.7] in the case where $a = 1$ in that
proof.
\qed\enddemo

\proclaim{Definition 1.1.4 ([\xgl, Definition 2.3])}
For $\l \in \Lambda$, define $\phi_{\l} : W(\l) \times W(\l)
\rightarrow R$ by $\phi_{\l}(C_S, C_T) = \phi_1(S, T)$, where $S, T
\in M(\l)$, and extend $\phi_{\l}$ bilinearly.
\endproclaim

From now on, we assume that $R$ is a field.
We now state one of the key results of the
paper [\xgl], which is the main one we aim to generalise to infinite
dimensions.  This result genuinely requires $\Lambda$ to be finite.

\proclaim{Theorem 1.1.5 (Graham, Lehrer)}
For $\l \in \Lambda$, the subspace $\rad(\l)$ of $W(\l)$ given by $$
\rad(\l) := \{x \in W(\l) : \phi_{\l}(x, y) = 0 \text{\rm \ for all \ } y
\in W(\l) \}
$$ is an $A$-submodule of $W(\l)$.  The quotient module $L(\l) :=
W(\l)/ \rad(\l)$ is absolutely irreducible.  

Define $\Lambda_0 := \{\l \in \Lambda : \phi_{\l} \ne 0\}$.  Then
the set $\{L(\l) : \l \in \Lambda_0 \}$ is a complete set
of equivalences classes of absolutely irreducible $A$-modules.
\endproclaim

\demo{Proof}
This follows from [\xgl, Proposition 3.2, Theorem 3.4].
\qed\enddemo

\vskip 20pt

\subhead 1.2 Obstructions to cellular algebras in
infinite dimensions \endsubhead

Note that in the definition of a cellular algebra, there was no
requirement that the set $\Lambda$ be finite, although the sets
$M(\l)$ for a fixed $\l$ must be finite.  (Removing the latter
assumption would mean that the cell modules would be infinite
dimensional, which would be very inconvenient.)  

The following simple example (which is also discussed in the
remarks following [\xkx, Theorem 3.1])
shows that Theorem 1.1.5 fails fairly badly
for a general cellular algebra with an infinite poset $\Lambda$.
\eject
\proclaim{Proposition 1.2.1}
Let $R$ be an algebraically closed field and let $A$ be the ring of
polynomials $R[x]$.  Let $\Lambda$ be the set of natural numbers
(including $0$), ordered by the reverse of the usual order.
Let $M(\l)$ be a set containing one element for
each $\l \in \Lambda$.  If $S \in M(\l)$, we define $C(S, S) =
x^{\l}$.  Let $*$ be the identity map on $A$.

Then $(\Lambda, M, C, *)$ is a cell datum for $A$ over $R$, in the
sense of Definition 1.1.1.
\endproclaim

\demo{Proof}
We omit the proof, since it is almost trivial.
\qed\enddemo

Although this algebra $R[x]$ satisfies the cellular axioms, Theorem
1.1.5 fails severely when applied to it.

\proclaim{Proposition 1.2.2}
The irreducible representations of $R[x]$ over $R$ are obtained by
evaluation at $x$ and are therefore in bijection with the elements of
$R$.

The set $\Lambda_0$ consists only of the element $0 \in \Lambda$.

Therefore Theorem 1.1.5 is false in general if $\Lambda$ is infinite.
\endproclaim

\demo{Proof}
The first assertion follows from Schur's lemma, since $R$ is
algebraically closed.

The second assertion follows from the fact that $x^r . x^r \equiv 0 \mod
\langle x^r \rangle$, unless $r = 0$ in which case 
$x^r . x^r \equiv x^r \mod \langle x^r \rangle$.

It is now obvious that the set $\Lambda_0$ does not index the
absolutely irreducible $R[x]$-modules over $R$, which is the third assertion.
\qed\enddemo

We observe that the only irreducible representation of $R[x]$ which
extends to a ``continuous''
representation of the power series ring $R[[x]]$ is the
one which is indexed by the element $0 \in \Lambda_0$.  This suggests
that in order to find an analogue of Theorem 1.1.5 in the infinite
dimensional case, one should consider modules for a suitable
completion of the original algebra.

\vskip 20pt

\head 2. Procellular algebras \endhead

From now on, we assume that $R$ is a field unless otherwise stated.

\subhead 2.1 Construction of procellular algebras \endsubhead

\proclaim{Definition 2.1.1}
A coideal $I$ of the poset $\Lambda$ is a subset of $\Lambda$ for which
$a \in I$ and $a < b$ imply $b \in I$.

The coideal $\langle x_1, x_2, \ldots, x_r \rangle$ 
generated by the elements $x_1, x_2, \ldots, x_r$ consists of the
elements of $\Lambda$ given by $$
\{ \l \in \Lambda: x_i \leq \l \text{\rm \ for some } 1 \leq i \leq r\}
.$$
\endproclaim

\proclaim{Definition 2.1.2}
We say that the cell datum $(\Lambda, M, C, *)$ for an algebra $A$
is of profinite type if $\Lambda$ is infinite and if for each $a \in
\Lambda$, the coideal $\langle a \rangle$ is a finite set.
\endproclaim

The coideals of the poset $\Lambda$ will be used to define cellular quotients of
the infinite dimensional algebra $A$.  We will then make an inverse
system using these quotients.  To achieve this, we make the following
definitions. 

\proclaim{Definition 2.1.3}
Let $\Lambda$ be the poset of a cell datum of profinite type.  We
denote the set of finite coideals of $\Lambda$, ordered by
inclusion, by $\Pi$.

If $P \in \Pi$, we write $A_P$ for the free $R$-module with basis
parametrised by the set $$
\{ C_P(S, T) : S, T \in M(\l), \l \in P \}
.$$  We write $I_P$ for the $R$-submodule of $A$ spanned by
all elements $C(S, T)$ where $S, T \in M(\l)$ for $\l
\not\in P$.  
\endproclaim
\eject
\proclaim{Proposition 2.1.4}
Let $P \in \Pi$ as above.  Then $A_P$ is a homomorphic image $\psi_P(A)$
of $A$ corresponding to the quotient of $A$ by the
ideal $I_P$.  The algebra $A_P$ is a finite dimensional cellular
algebra which inherits a cell datum from $A$ by
restriction of the poset $\Lambda$ to the set $P$.
\endproclaim

\demo{Proof}
This follows from the defining axioms for a cellular algebra.
\qed\enddemo

\proclaim{Definition 2.1.5}
Maintain the above notation.  Let $P_1, P_2 \in \Pi$ and suppose $P_1
\supset P_2$.  We define the $R$-linear map $\psi_{P_1, P_2}: A_{P_1}
\rightarrow A_{P_2}$ via the conditions $$
\psi_{P_1, P_2}(C_{P_1}(S, T)) = \cases
C_{P_2}(S, T) & \text{ if } S, T \in M(\lambda) \text{ where } 
\lambda \in P_2,\cr
0 & \text{ otherwise.}\cr
\endcases$$
\endproclaim

\proclaim{Lemma 2.1.6}
The maps $\psi$ of Definition 2.1.5 are surjective algebra homomorphisms.
\endproclaim

\demo{Proof}
This follows from the fact that $\psi_{P_1, P_2} . \psi_{P_1} =
\psi_{P_2}$.
\qed\enddemo

We now use the maps $\psi$ to define a inverse limit of
cellular algebras.  This idea has its roots in [\xthesis, \S6.4].

\proclaim{Definition 2.1.7}
Let $A$ be a cellular algebra with a cell datum of profinite type.
Consider the inverse system whose
elements are the sets $A_P$ for $P \in \Pi$ and whose
homomorphisms are given by $$
\psi_{P_1, P_2} : A_{P_1} \rightarrow A_{P_2}
$$ whenever
$P_1 \supset P_2$.  Denote the associated inverse limit by $\widehat
A$.  We call an algebra $\widehat A$ arising in this way a procellular
algebra.
\endproclaim

\demo{Remark}
The terminology ``procellular'' is by analogy with profinite groups,
which are inverse limits of finite groups.
\enddemo

\vskip 20pt

\subhead 2.2 Properties of procellular algebras \endsubhead

An element of a procellular algebra $\widehat A$
can be regarded as a certain infinite combination of
the basis elements $C(S, T)$ of $A$, as follows.

\proclaim{Proposition 2.2.1}
There is a canonical bijection between elements of $\widehat A$ and
formal infinite sums $$
\sum_{(T, T')} a(T, T') C(T, T')
,$$ as $(T, T')$ ranges over the domain of $C$.  

The element of
$\widehat A$ corresponding to this sum projects in the inverse system
to the element of $A_P$ given by $$
\sum_{(T, T')} a(T, T') \, \psi_P(C(T, T'))
.$$
\endproclaim

\demo{Remark}
The properties of the map $\psi_P$ (see Proposition 2.1.4) mean that 
the second sum appearing above is finite.
\enddemo

\demo{Proof}
By arguments familiar from \S2.1, we see that the elements of $A_P$ as
defined above give (as $P$ varies over $\Pi$) an element of $\widehat
A$.  It is clear that this correspondence is an injection.

It remains to show that every element of $\widehat A$ is of this form.
Consider an element $C(T, T') \in A$, where $T, T' \in M(\l)$.
Note that $\langle \l \rangle \in \Pi$.  Define $a(T, T')$ to be the
coefficient of $\psi_{\langle \l \rangle}
(C(T, T'))$ in $A_{\langle \l \rangle}$.  Then for any
$P$ such that $\l \in P$, the coefficient of $\psi_P(C(T, T'))$ in $A_P$ must
also be $a(T, T')$, by consideration of $\psi_{P, \langle \l
\rangle}$.  It follows that $\widehat A$ is determined solely by the
coefficients $a(T, T')$, which proves the assertion.
\qed\enddemo

\proclaim{Lemma 2.2.2}
The procellular algebra $\widehat A$ admits an algebra anti-auto\-morph\-ism
$\widehat *$ of order 2 which arises from the anti-auto\-morph\-isms $*_P$ on
the cellular algebras $A_P$.
\endproclaim

\demo{Proof}
This is immediate from the construction of the inverse system.
\qed\enddemo

\proclaim{Definition 2.2.3}
Let $P \in \Pi$.  We define $\widehat I_P$ to be the ideal of
$\widehat A$ given by the infinite sums $$
\sum a(T, T') C(T, T')
$$ where $a(T, T') \ne 0 \Rightarrow T, T' \in M(\l)$ for $\l \not\in P$.
\endproclaim

The procellular algebra $\widehat A$ is equipped with a natural
topology arising from the poset $\Lambda$.  It will turn out that
$\widehat A$ is in fact the completion of $A$ with respect to this topology.

\proclaim{Proposition 2.2.4}
Let the set $\{\widehat I_P : P \in \Pi\}$ be a base of neighbourhoods of $0
\in \widehat A$.  This gives $\widehat A$ the structure of a Hausdorff,
complete topological ring in which the operation $\widehat *$ is a 
homeomorphism.
\endproclaim

\demo{Remark}
By a topological ring, we mean a ring in which addition, negation and
multiplication are continuous operations, so that in particular,
addition is a homeomorphism.  The topology may be specified by giving
a base of neighbourhoods for $0$, meaning that any open neighbourhood
of $0$ contains one of the neighbourhoods in the base.  The open sets
are those arising from the structure of the ring as a topological
abelian group under addition.
\enddemo

\demo{Proof}
The assertion that $\widehat A$ acquires the structure of a
topological ring is clear, because the neighbourhoods in the base are ideals.

We show that the intersection of the ideals $\widehat I_P$ is
trivial.  Let $0 \ne x \in \bigcap_{P \in \Pi} \widehat I_P$, and 
suppose $C(T, T')$
appears with nonzero coefficient $a(T, T')$ in the infinite sum
expansion of $x$.  If $T, T' \in M(\l)$, then $\psi_{\langle \l
\rangle}(x) \ne 0$, and therefore $x \not\in \widehat I_{\langle \l
\rangle}$, which is a contradiction.

It follows that $\widehat I_{\langle \l \rangle}$ and $x + \widehat 
I_{\langle \l
\rangle}$ are disjoint open sets which contain the points $0$ and $x$
respectively.  We deduce that the topology is Hausdorff.
Completeness can be checked from the infinite sum realisation of $\widehat A$.

Finally we observe that $\widehat *$ is of order 2, is an algebra
homomorphism, and fixes the ideals $\widehat I_P$ setwise.  This
proves the last assertion.
\qed\enddemo

\proclaim{Corollary 2.2.5}
The algebra $A$ embeds canonically in $\widehat A$ by considering the
finite sums of $C(T, T')$ as special cases of the infinite sums of
Proposition 2.2.1.  Furthermore, $\widehat A$ is the completion of $A$
with respect to the topology defined above.
\endproclaim

\demo{Proof}
It is clear that the embedding works as stated.  We also observe that
any infinite sum can be approximated arbitrarily closely by a finite
sum in the given topology, because all sets $P \in \Pi$ are finite.
It follows that $A$ embeds densely and thus that $\widehat A$ is the
completion of $A$.
\qed\enddemo

\vskip 20pt

\subhead 2.3 Representations of procellular algebras \endsubhead

As might be expected from the theory of topological groups, we will
concentrate on the ``smooth'' representations of the procellular
algebras.  We are also only interested in finite dimensional
representations.  We define the smooth representations as follows,
based on [\xcarti, Definition 1.1].

\proclaim{Definition 2.3.1}
A representation $\rho$ of $\widehat A$ on a space $V$ 
is smooth if the annihilator of every vector of $V$ is open.
\endproclaim

\proclaim{Lemma 2.3.2}
Let $\rho$ be a finite dimensional representation of $\widehat A$ on a
space $V$.  Then $\rho$ is smooth if and only if $\ker \rho$ is open.
\endproclaim

\demo{Proof}
Suppose $\rho$ is smooth.
Then $\ker \rho$ is the intersection of a finite number of
annihilators of vectors $v$ and is therefore open.

Conversely, assume $\ker \rho$ is open, and
pick $v \in V$.  Then the annihilator of $v$ is a union of cosets of
$\ker \rho$, and is therefore open.
\qed\enddemo

We will also refer to ``smooth modules'', with the obvious meaning.

\proclaim{Definition 2.3.3}
Let $\l \in \Lambda$.  The cell module $\widehat W(\l)$ for $\widehat
A$ is the left module obtained by considering the cell module
$W(\l)$ for $A_{\langle \l \rangle}$ as a module for $\widehat A$ via
the homomorphism $\psi_{\langle \l \rangle}$.

The module $\widehat W(\l)$ inherits the bilinear form $\phi_{\l}$
from $W(\l)$.
We define $\widehat L(\l)$ as the $\widehat A$-module
corresponding to the $A_{\langle \l \rangle}$-module $$
L(\l) = W(\l)/ \rad(\l).$$
\endproclaim

The importance of the modules $\widehat L(\l)$ is demonstrated by the
following theorem (compare with Theorem 1.1.5).

\proclaim{Theorem 2.3.4}
Let $A$ be a cellular algebra with a datum of profinite type, and let
$\widehat A$ be the corresponding procellular algebra.  

Define $\Lambda_0 := \{\l \in \Lambda : \phi_{\l} \ne 0\}$. 

Then the set $\{\widehat L(\l) : \l \in \Lambda_0 \}$
is a complete set
of equivalences classes of absolutely irreducible smooth $\widehat A$-modules.
\endproclaim

\demo{Proof}
It is clear from the properties of $L(\l)$ that the module $\widehat
L(\l)$ is absolutely irreducible for each $\l$.

To check the smoothness property, we observe that $\widehat I_{\langle
\l \rangle}$ is an open set which lies in the kernel of the
representation corresponding to $\widehat W(\l)$.  
The kernel of $\widehat L(\l)$ is a union of cosets of this open
set and is therefore open.  Thus by Lemma 2.3.2, the module is smooth.

We see that two modules $\widehat L(\l)$ and $\widehat L(\mu)$ are
nonisomorphic by considering them as modules for the finite
dimensional cellular algebra $A_{\langle \l, \mu \rangle}$.

It remains to prove that the given modules exhaust the set of
absolutely irreducible smooth $\widehat A$-modules.  Let $L$ be such a
module.  Then the kernel of $L$ is open, so $\widehat I_P$ annihilates
$L$ for some $P \in \Pi$.  Thus $L$ is an absolutely
irreducible $A_P$-module in a canonical way.  
Proposition 2.1.4 shows that $A_P$ is naturally a finite dimensional 
cellular algebra,
which means by Theorem 1.1.5 that $L$ is a cell module for $A_P$
corresponding to some $\mu \in P$.  The definition of the modules
$\widehat L(\l)$ means that $L$, regarded as an $\widehat A$-module,
is nothing other than $\widehat L(\mu)$.
\qed\enddemo

\vskip 20pt

\head 3. Examples \endhead

It is now clear that the power series counterexample in Proposition
1.2.2 is compatible with Theorem 2.3.4.  Any smooth module $L$ for
the procellular algebra $\widehat{R[x]} = R[[x]]$ 
must satisfy the property that $x^n . L = 0$ for
sufficiently large $n$, meaning that $x . L = 0$.  Thus the only
smooth module for the procellular algebra corresponds to the element
$0 \in \Lambda_0$.

This is a trivial example, but there are also interesting examples of
procellular algebras as we explain in the next section.

\subhead 3.1 The algebras $\dot U$ \endsubhead

In [\xlusc, \S29], Lusztig essentially proves that the algebras $\dot
U$, which are associated with Dynkin diagrams of various types, are
cellular.  (This result is referred to as the Peter--Weyl theorem.)  
The analogues of the finite dimensional algebras $A_P$ are also considered.  
The algebras $\dot U$ are equipped with canonical bases, $\dot B$, which
have many beautiful properties which are developed in [\xlusc].

We concentrate in this
section on the case where $\dot U$ corresponds to a Dynkin diagram of
type $A$, because this has interesting connections with the $q$-Schur
algebra, and because it is better understood.  Since we are only
introducing these algebras for illustrative purposes, we shall not go
into the full details.  The interested reader is referred to the
literature for more information.

The $q$-Schur algebra $S_q(n, r)$, which first appeared in [\xdjb], is
a finite dimensional algebra over a base ring containing an
indeterminate $q$ and its square root $q^{1/2}$.  This algebra was
proved by the author to be cellular with respect to a Murphy-type basis 
[\xthesis, Proposition 6.2.1].  Du [\xdub] introduced a 
canonical basis $\{\theta_{S, T}\}$ for this algebra, and the $q$-Schur
algebra is also cellular with respect to the $\theta$-basis.

\proclaim{Proposition 3.1.1}
The $q$-Schur algebra $S_q(n, r)$ has a basis $\{\theta_{S, T}\}$,
where $S, T$ are semistandard tableaux of the same shape with $r$
boxes each and entries from the integers $\{1, 2, \ldots, n\}$.

This gives the $q$-Schur algebra the structure of a cellular algebra,
where $\Lambda = \Lambda^+(n, r)$ ordered by the dominance order (as
in [\xjagb, \S2]), $M(\l)$ is the set of semistandard tableaux of shape
$\l$ (as above), $C$ is the map $C(S, T) = \theta_{S, T}$ and $*$ is
the anti-auto\-morph\-ism sending $\theta_{S, T}$ to $\theta_{T, S}$.
\endproclaim

\demo{Remark}
Recall that ``semistandard'' means that the entries increase weakly
along the rows, and strictly down the columns of the tableau.  

This result is a generalisation of the example of the Hecke algebra of
type $A$ appearing in [\xgl, Example 1.2].
\enddemo

\demo{Proof}
Apart from the assertion involving $*$, this is simply [\xdr, Theorem
5.3.3].  

It follows from the definition of the $\theta$-basis given in [\xdub] that
there is an anti-auto\-morph\-ism sending $\theta_{\l, \mu}^w$ to
$\theta_{\mu, \l}^{w^{-1}}$, where $\l, \mu$ are compositions of $r$
into at most $n$ pieces, and $w$ is a certain double coset
representative in the symmetric group $S_r$.

The assertion about $*$ comes from properties of the
Robinson--Schensted correspondence used in [\xdr, \S5.2].  It is a
well-known property of this correspondence that inversion of the group
element corresponds to exchange of the pair of tableaux.  In the
notation of semistandard tableaux, the anti-auto\-morph\-ism of the
previous paragraph is nothing other than $*$.
\qed\enddemo

The terminology ``canonical'' for this basis is justified, as the
following result shows.

\proclaim{Proposition 3.1.2 (Du)}
For each $r$, there is an epimorphism of algebras $$
\pi_r : \dot U(sl_n) \rightarrow S_q(n, r).
$$  This takes elements of the basis $\dot B$ to zero or
to canonical basis elements $\theta_{S, T}$.  

For each basis element $b \in \dot B$, there exists an $r$ such that
$\pi_r(b) \ne 0$.
\endproclaim

\demo{Proof}
This follows from [\xduf, Theorem 3.5] and the remarks in [\xduh, \S5.5].
\qed\enddemo

\proclaim{Proposition 3.1.3}
Let $\Lambda$ be the set of weights associated to finite dimensional 
simple modules for
the Lie algebra $sl_n$ over ${\Bbb C}$, ordered by the usual dominance
order, where $\l < \mu$ means $\l$ dominates $\mu$.  

Let $M(\l)$ be
the set of semistandard tableaux with $l_i$ boxes in the $i$-th row,
where $l_i \geq l_j$ for $i < j$, $l_n = 0$ and $\l_i = l_i -
l_{i+1}$.

Let $C$ be the map taking a pair $(S, T)$ of semistandard tableaux of
the same shape to the element $b \in \dot B$ corresponding to the
canonical basis element $\theta_{S, T}$ of the $q$-Schur algebra $S_q(n,
r)$ where $r = \sum \l_i$.

Let $*$ be the map induced by sending $\theta_{S, T}$ to $\theta_{T,
S}$.

Then $(\Lambda, M, C, *)$ is a cell datum for $\dot U$, and the datum
is of profinite type.
\endproclaim

\demo{Note}
One can check this directly using the theory of quantum groups
developed in [\xlusc], but the proof we sketch below is based on
inverse systems of $q$-Schur algebras.
\enddemo

\demo{Proof}
The cell structure is essentially inherited from the cellular
structure of the algebras $S_q(n, r)$ given in Proposition 3.1.1 and
the homomorphisms of Proposition 3.1.2.  

As mentioned in [\xduh, \S5.5], one can
construct an inverse limit of $q$-Schur algebras using
the epimorphisms $\psi_{i + jn}$ 
from $S_q(n, i + (j+1)n)$ to $S_q(n, i + jn)$ which
are dual to the $q$-analogue of coalgebra injection corresponding to
multiplication by the determinant map.  The effect of $\psi$
on the $\theta$-basis is known explicitly [\xduh, 5.4
(b)] and takes a $\theta$-basis element to another $\theta$-basis
element or to zero.  The properties of cell ideals in cellular
algebras now imply that if $\theta_{S_1, T_1}$ and $\theta_{S_2, T_2}$
are in the same cell (\idest $S_1, T_1, S_2, T_2 \in M(\l)$ for the
same $\l$), then their images under $\psi$ are either both
zero, or both nonzero and in the same cell.

It now follows that $\dot U$ does inherit a cellular structure from
the inverse system described above.
We can thus identify a basis element $b$ with its images under the
maps $\psi$, the ``lowest'' of which corresponds to a basis element
$\theta$ whose tableaux have leftmost columns with
fewer than $n$ entries.  The fact
that the dominance order is the relevant order can be deduced from
properties of 
the poset order in $S_q(n, r)$.  Using these facts, we can check that
the assertions in the statement hold.
\qed\enddemo

\demo{Remark}
One can check from the Robinson--Schensted 
correspondence that the effect of $\psi$ is as follows.
If $S$ and $T$ have a column of length $n$ then $S', T'$ are obtained from
$S, T$ respectively by removal of the leftmost columns, and $\theta_{S,
T}$ maps to $\theta_{S', T'}$.
If this is not the case then $\theta_{S, T}$ is mapped to zero.

This phenomenon of column removal was also studied
in [\xthesis, \S6.4].
\enddemo

\vskip 20pt
\eject
\subhead 3.2 Further remarks on procellular algebras \endsubhead

The cell datum for $\dot U$ given in \S3.1 gives rise to a procellular
algebra $\widehat U$ which has interesting properties.  Du [\xduh,
\S5.5] calls the algebra $\widehat U$ the completion of $U$ with
respect to $q$-Schur algebras.  The algebra $\widehat U$ has some
convenient properties not shared by $\dot U$, such as the existence of
a multiplicative identity.

The algebra $\dot U$ has a positivity property for its structure constants
with respect to $\dot B$, which can be easily seen (by the results in
\S3.1) to be equivalent to the positivity of the structure constants
for the $\theta$-basis.  The positivity property for the
$\theta$-basis was proved by the author in [\xrmge].
In both cases, the structure constants have interpretations in the geometrical
framework of perverse sheaves.

The completion $\widehat U(sl_n)$, as well as containing the algebra 
$\dot U(sl_n)$, contains the quantized enveloping
algebras $U(sl_n)$ in the sense of Drinfel'd and Jimbo.  This was
shown in [\xthesis, Theorem 6.4.19].

Several more of the main results concerning cellular algebras go over
to procellular algebras essentially verbatim, so we will not discuss
them explicitly here.

Another interesting example of a procellular algebra is the
completion of the affine Temperley--Lieb algebra.  In this case,
the smooth modules for the completed algebra correspond to finite
dimensional indecomposable modules of the original algebra.
Furthermore, the procellular completion depends on a nonzero
parameter; altering this parameter produces different families of 
indecomposable modules.  The reader is referred to [\xeg] for
further details.

\vskip 20pt

\head Acknowledgements \endhead

The author is grateful to Steffen K\"onig and Changchang Xi for 
conversations and correspondence relevant to this work, and to the
referee for helpful suggestions.

\vskip 1cm
\head References \endhead

\item{[\xcarti]}
{P. Cartier, {\it Representations of $p$-adic groups: a survey},
Proc. Sympo. Pure Math., {\bf 33} (1979), part 1, 111-155}
\item{[\xdjb]}
{R. Dipper and G.D. James, {\it The $q$-Schur algebra}, Proc. L.M.S. {\bf 59}
(1989), 23--50.}
\item{[\xdub]}
{J. Du, {\it Kazhdan-Lusztig bases and isomorphism theorems for $q$-Schur
algebras}, Contemp. Math. {\bf 139} (1992), 121--140.}
\item{[\xduh]}
{J. Du, {\it $q$-Schur Algebras, Asymptotic Forms, and Quantum $SL_n$},
J. Alg {\bf 177} (1995), 385--408.}
\item{[\xduf]}
{J. Du, {\it Global IC Bases for Quantum Linear Groups}, J. Pure
Appl. Alg. {\bf 114} (1996), 25--37.}
\item{[\xdr]}
{J. Du and H. Rui, {\it Based algebras and standard bases for quasi-hereditary
algebras}, preprint.}
\item{[\xeg]}
{K. Erdmann and R.M. Green, {\it On Representations of Affine
Temperley--Lieb Algebras, II}, Pac. J. Math., to appear}
\item{[\xgl]}
{J.J. Graham and G.I. Lehrer, {\it Cellular Algebras}, Invent. Math. {\bf 123}
(1996), 1--34.}
\item{[\xjagb]}
{J.A. Green,
{\it Combinatorics and the Schur algebra}, J. Pure Appl. Alg. {\bf 88} (1993)
89--106.}
\item{[\xthesis]}
{R.M. Green, Ph.D. thesis, University of Warwick, 1995.}
\item{[\xrmge]}
{R.M. Green, {\it Positivity Properties for $q$-Schur Algebras},
Proc. Cam. Phil. Soc., {\bf 122} (1997), 401--414.}
\item{[\xkx]}
{S. K\"onig and C. Xi, 
{\it On the structure of cellular algebras}, to appear in the 
Proceedings of the 8th International Conference on Representations of
     Algebras.}
\item{[\xlusc]}
{G. Lusztig, {\it Introduction to Quantum Groups}, Birkh\"auser, 1993.}

\end